\documentclass[11pt]{article}
\usepackage{amssymb,amsmath,amsfonts}
\usepackage{graphics}
\usepackage[mathscr]{eucal}

\allowdisplaybreaks  

\newcommand{\sinush}{\mathop{\rm sinh}}
\newcommand{\cosinh}{\mathop{\rm cosh}}
\newcommand{\tgh}{\mathop{\rm tanh}}
\newcommand{\arctanh}{\mathop{\rm artanh}}
\newcommand{\arch}{\mathop{\rm arcosh}}

\begin{document}

\title{\Large{AN INFINITESIMALLY NONRIGID POLYHEDRON WITH NONSTATIONARY VOLUME IN THE LOBACHEVSKY 3-SPACE}}

\author{Dmitriy~Slutskiy \thanks{The author is supported in part by
the Council of Grants from the President of the Russian Federation
(Grant NSh-6613.2010.1), by the Federal Targeted Programme on
Scientific and Pedagogical-Scientific Staff of Innovative Russia for
2009-2013 (State Contract No. 02.740.11.0457) and by the Russian
Foundation for Basic Research (Grant 10-01-91000-ANF\_a).}}

\date{}

\maketitle


\begin{abstract}
We give an example of an infinitesimally nonrigid polyhedron in the
Lobachevsky 3-space and construct an infinitesimal flex of that
polyhedron such that the volume of the polyhedron isn't stationary
under the flex.
\end{abstract}

\begin{center}

\end{center}

\section{Introduction}

The Bellows Conjecture states that every flexible polyhedron
preserves its oriented volume during the flex. In 1996
I.\,Kh.~Sabitov \cite{Sabitov1996} gave an affirmative answer to the
Bellows Conjecture in the Euclidean 3-space. In 1997
V.\,A.~Alexandrov \cite{Alexandrov1997} has built a flexible
polyhedron in the spherical 3-space which changes its volume during
the flex. The question whether the Bellows Conjecture holds true in
the Lobachevsky 3-space is still open.

In the note of the editor of the Russian translation of
\cite{Connelly1980} I.\,Kh.~Sabitov proposed to consider the Bellows
Conjecture at the level of infinitesimal flexes. Roughly, we can
formulate I.\,Kh.~Sabitov's question as follows: is it true that,
for every infinitesimally nonrigid polyhedron, the volume it bounds
is stationary under its infinitesimal flex? In case the answer to
I.\,Kh.~Sabitov's question were positive, we would automatically
validate the Bellows Conjecture for the flexible polyhedra. Of
course, we can always additionally triangulate any initial face of a
polyhedron so that there exists a new vertex $A$ of the
triangulation which is an internal point of the initial face, then
attach to $A$ a nonzero velocity vector orthogonal to the initial
face, leave all other vertices of the polyhedron fixed and thus
construct an infinitesimal flex of the new polyhedron based on the
movements of all its vertices. The volume of the polyhedron with the
``false'' vertex under the constructed infinitesimal flex is
nonstationary, but this trivial example is of a little interest to
study the Bellows Conjecture.

Having constructed a nontrivial counterexample in
\cite{Alexandrov1989}, V.\,A.~Alexandrov gave a negative answer to
I.\,Kh.~Sabitov's question for infinitesimally nonrigid polyhedra in
the Euclidean 3-space. An example of a flexible polyhedron in the
spherical 3-space, constructed in \cite{Alexandrov1997}, which
changes its volume during the flex, yields that the answer to this
question is also negative for infinitesimally nonrigid polyhedra in
the spherical 3-space. The main result of this paper reads as
follows.

\textbf{Theorem} \emph{In the Lobachevsky 3-space there is a
sphere-homeomorphic intersection-free polyhedron and its
infinitesimal flex such that the volume it bounds isn't stationary
under the flex.}

The polyhedron mentioned in the theorem is built explicitely. It's
similar to a polyhedron in the Euclidean 3-space which was first
constructed by A.\,D.~Alexandrov and S.\,M.~Vladimirova
\cite{Alexandrov_Vladimirova1962} and later studied by A.\,D. Milka
\cite{Milka2000}.

\section{Constructing $\mathscr{S}$} \label{paragraph_polyhedron_construction}

Throughout this paper we call a polyhedral surface a polyhedron.

\begin{figure}
\includegraphics{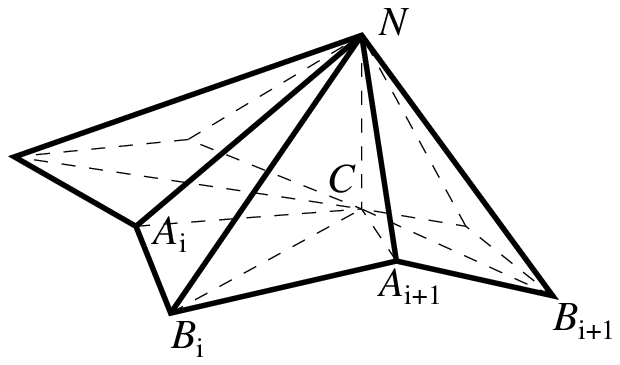} \hfill
\includegraphics{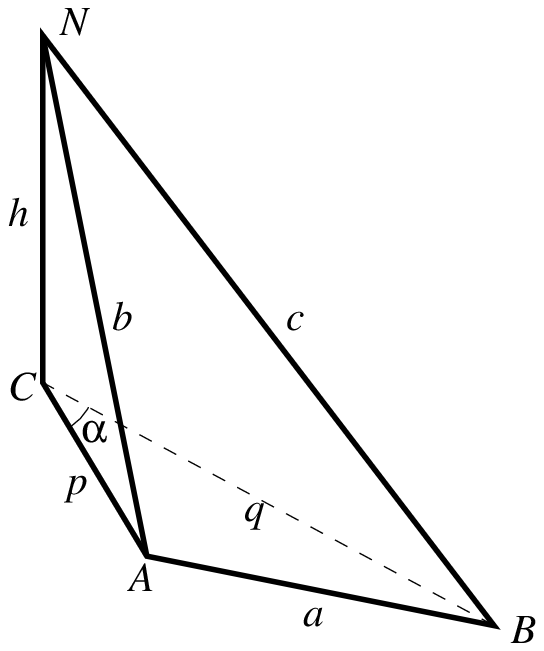} \\
\parbox[t]{7cm}{ \caption{The lateral surface of $\mathscr{P}$.}\label{sss_whole}} \hfill
\parbox[t]{7cm}{\caption{The tetrahedron $\mathscr{T}$.}\label{sss_part}}
\end{figure}

Consider a regular pyramid $\mathscr{P}$ in the Lobachevsky 3-space
with a regular concave star with $n$ petals as the base. We denote
vertices of the star by $A_i$, $B_i$, $i = 1,...,n$, and we note
that the orthogonal projection of the vertex $N$ of $\mathscr{P}$
onto its base coincides with the center $C$ of the star, see
Fig.~\ref{sss_whole}. We reflect $\mathscr{P}$ in the plane that
contains its base and denote by $\mathscr{S}$ a suspension which
consists of both initial and reflected pyramids without their common
base. We denote by $S$ the vertex of $\mathscr{S}$ symmetric to $N$
with respect to the plane containing the base of $\mathscr{P}$. A
cycle formed by the edges of the base of $\mathscr{P}$ is called the
equator of the suspension $\mathscr{S}$.

Note that the lengths of all edges of the equator of $\mathscr{S}$
are equal to each other by construction. Moreover, the lengths of
all edges $SA_i$, $NA_i$, $i = 1,...,n$, are equal to each other,
and also the lengths of all edges $NB_i$, $SB_i$, $i = 1,...,n$, are
equal to each other too.

By construction, $\mathscr{S}$ possesses multiple symmetries and the
spatial body bounded by $\mathscr{S}$ consists of identical
tetrahedral ``bricks''. Consider one of these tetrahedra, see
Fig.~\ref{sss_part}. Denote its surface by $\mathscr{T}$, and its
vertices by $N$, $A$, $B$, $C$. Note that $\angle ACN = \angle BCN =
\pi/2$ by construction. Let's use the following notations for the
lengths of the edges and for the plane angles of $\mathscr{T}$:
$|CN| = h$, $|CA| = p$, $|CB| = q$, $|AB| = a$, $|NA| = b$, $|NB| =
c$, $\angle ACB = \alpha$, $\angle CAN = \beta$, $\angle BAN =
\gamma$, $\angle CAB = \delta$, $\angle CBN = \varphi$, $\angle CBA
= \psi$, $\angle ABN = \theta$, $\angle ANB = \lambda$, $\angle CNA
= \mu$, $\angle CNB = \nu$. Denote the dihedral angles of
$\mathscr{T}$ at the edge $AB$ by $\angle AB$, at the edge $NA$ by
$\angle NA$, and at the edge $NB$ by $\angle NB$.

By construction, the dihedral angle of $\mathscr{T}$ at the edge
$CN$ is equal to $\alpha$, the dihedral angles of $\mathscr{S}$ at
the edges of its equator are equal to $2 \angle AB$, at the edges
$NA_i$ and $SA_i$ , $i = 1,...,n,$ are equal to $2 \angle NA$, and
at the edges $NB_i$ and $SB_i$ , $i = 1,...,n,$ are equal to $2
\angle NB$.

Further we show that the suspension $\mathscr{S}$ constructed above
can be taken as a polyhedron whose existence is proclaimed by the
theorem.

\section{A condition for infinitesimal nonrigidity} \label{paragraph_flex_condition}

A deformation of a polyhedral surface $\mathscr{S}$ is a family of
surfaces $\mathscr{S}(t)$, $t \in (-1,1),$ which depends
analytically on the parameter $t$, preserves the combinatorial
structure of $\mathscr{S}$, and is such that $\mathscr{S}(0) =
\mathscr{S}$.

A deformation of a polyhedral surface $\mathscr{S}$ with triangular
faces is called its infinitesimal flex if the lengths of all edges
of $\mathscr{S}(t)$ are stationary at $t = 0$.

An infinitesimal flex is called nontrivial if there exist two
vertices of $\mathscr{S}(t)$ which are not connected by an edge of
$\mathscr{S}(t)$ and are such that the spatial distance between them
is not stationary.

A polyhedron is called infinitesimally nonrigid if it possesses a
nontrivial infinitesimal flex.

Determine a deformation of the suspension $\mathscr{S}$ constructed
in the previous section as follows. The point $C$ is fixed. At the
moment $t$, the point $N$ goes to the point $N(t)$ lying on the ray
$\overrightarrow{CN}$ at the distance from $C$ determined by the
formula
\begin{equation}
\label{CN(t)} h(t) = h + tu,
\end{equation}
where $u$ is a real number which has a meaning of velocity and which
will be specified below. The point $S$ goes to the point $S(t)$
lying on the ray $\overrightarrow{CS}$ at the distance from $C$
determined by the formula (\ref{CN(t)}). The point $A_i$, $i =
1,...,n,$ goes to the point $A_i(t)$ lying on the ray
$\overrightarrow{CA_i}$ at the distance from $C$ determined by the
formula $ \label{CA(t)} p(t) = p + tv$, where $v$ is a real number
which has a meaning of velocity. The point $B_i$, $i = 1,...,n,$
goes to the point $B_i(t)$ lying on the ray $\overrightarrow{CB_i}$
at the distance from $C$ determined by the formula $ \label{CB(t)}
q(t) = q + tw$, where $w$ is a real number which has a meaning of
velocity and which will be specified below.

In order to determine the movements of other points of the
suspension $\mathscr{S}(t)$ let's use the statement of Ceva's
theorem in the Lobachevsky space \cite{Prasolov2004}:

\emph{Given a triangle} $\triangle ABC$ \emph{and points} $\tilde{A}
$, $\tilde{B}$, \emph{and} $\tilde{C}$ \emph{that lie on sides}
$BC$, $CA$, \emph{and} $AB$ \emph{of} $\triangle ABC$. \emph{Then
the segments} $A\tilde{A}$, $B\tilde{B}$, \emph{and} $C\tilde{C}$
\emph{intersect at one point if and only if one of the following
equivalent relations holds}:
\begin{equation*} \label{ceva_angles}
\frac{\sin \angle AC\tilde{C}}{\sin \angle \tilde{C}CB} \frac{\sin
\angle BA\tilde{A}}{\sin \angle \tilde{A}AC} \frac{\sin \angle
CB\tilde{B}}{\sin \angle \tilde{B}BA}= 1;
\end{equation*}
\begin{equation} \label{ceva_lengths}
\frac{\sinush A\tilde{C}}{\sinush \tilde{C}B} \frac{\sinush
B\tilde{A}}{\sinush \tilde{A}C} \frac{\sinush C\tilde{B}}{\sinush
\tilde{B}A} = 1.
\end{equation}

In terms of the statement of Ceva's theorem, let's take the point
$P(t)$ of the segment $A(t)B(t)$ for which the equality
\begin{equation*} \label{point_on edge}
\frac{\sinush A(t)P(t)}{\sinush P(t)B(t)} = \frac{\sinush
AP}{\sinush PB}
\end{equation*}
holds true, as a new position of any point $P$ of the edge $AB$ at
the moment $t$.

To determine the movement of an internal point $Q$ of the face
$\triangle ABC$, at first we construct points $\tilde{A} $,
$\tilde{B}$, and $\tilde{C}$, as the intersections of the edges
$BC$, $CA$, and $AB$ with the rays $AQ$, $BQ$, and $CQ$, and then
determine their positions $\tilde{A}(t)$, $\tilde{B}(t)$, and
$\tilde{C}(t)$ at the moment $t$ by the method described above. By
Ceva's theorem, the segments $A(t)\tilde{A}(t)$, $B(t)\tilde{B}(t)$,
and $C(t)\tilde{C}(t)$ intersect at one point (the relation
(\ref{ceva_lengths}) remains true at every moment $t$). Consider
this point of intersection as a new position $Q(t)$ of the point $Q$
at the moment $t$.

The deformation of $\mathscr{S}$ described above, naturally produces
a deformation of the tetrahedron $\mathscr{T}$ which we denote by
$\mathscr{T}(t)$. The lengths of all edges as well as the values of
all plane and dihedral angles of $\mathscr{T}$ are functions in $t$
and their notations naturally succeed from the notations for the
corresponding entities of $\mathscr{T}$. For example, we denote the
length of the edge $N(t)A(t)$ by $b(t)$, the value of the plane
angle $\angle CA(t)N(t)$ by $\beta(t)$, and the value of the
dihedral angle of $\mathscr{T}(t)$ at the edge $N(t)A(t)$ by $\angle
N(t)A(t)$, etc.

Let's find a relation between $u$, $v$, and $w$ implying that the
deformation $\mathscr{S}(t)$ is an infinitesimal flex. We only need
to study the deformation of the face $ABN$ in $\mathscr{T}$ because
all faces of $\mathscr{S}$ move in the same way.

Apply the Pythagorean theorem for the Lobachevsky space
\cite{Vinberg1988} to the triangle $\triangle N(t)CA(t)$:
\begin{equation}
\label{cosinh_b(t)} \cosinh b(t) = \cosinh (h + tu) \cosinh (p + tv)
\end{equation}
and to the triangle $\triangle N(t)CB(t)$:
\begin{equation} \label{cosinh_c(t)}
\cosinh c(t) = \cosinh (h + tu) \cosinh (q + tw)
\end{equation}
of $\mathscr{T}(t)$.

Using the Cosine Law for the Lobachevsky space \cite{Vinberg1988}
applied to the triangle $\triangle A(t)CB(t)$, and taking it into
account that the angle $\alpha$ remains constant during the
deformation (and is equal to $\frac{\pi}{n}$), we get:
\begin{equation}
\label{cosinh_a(t)} \cosinh a(t) = \cosinh (p + tv) \cosinh (q + tw)
- \sinush (p + tv) \sinush (q + tw)\cos \alpha.
\end{equation}

Further it will be useful for us to study stationarity of the
function $f(t) = \cosinh l(t)$ instead of stationarity of the length
$l(t)$ of any edge of $\mathscr{S}(t)$, because $ \label{length_der}
f'(0) = l'(0) \sinush l(0) $ and $l(0) > 0$, and thus $f'(0) = 0$ if
and only if $l'(0) = 0$.

Let's differentiate (\ref{cosinh_b(t)}): $ \label{cosinh'_b(t)}
(\cosinh b(t))' = u \sinush (h + tu) \cosinh (p + tv) + v \cosinh (h
+ tu) \sinush (p + tv)$. Thus, stationarity of the length $b(t)$ of
the edge $N(t)A(t)$ is equivalent to the condition $
\label{cosinh'_b(0)} (\cosinh b(t))'|_{t=0} = u \sinush h \cosinh p
+ v \cosinh h \sinush p = 0$, or
\begin{equation}
\label{v(u)} v = - \frac{\tgh h}{\tgh p} u.
\end{equation}
Similarly, stationarity of the length $c(t)$ of the edge $N(t)B(t)$
is equivalent to the condition
\begin{equation}
\label{w(u)} w = - \frac{\tgh h}{\tgh q} u.
\end{equation}
Differentiating (\ref{cosinh_a(t)}), we find the condition for
stationarity of the length $a(t)$ of the edge $A(t)B(t)$:
\begin{equation}
\label{cosinh'_a(0)} (\cosinh a(t))'|_{t=0} = v \sinush p \cosinh q
+ w \cosinh p \sinush q - \cos \alpha  \{ v \cosinh p \sinush q + w
\sinush p \cosinh q \} = 0.
\end{equation}
Substituting (\ref{v(u)}) and (\ref{w(u)}) into
(\ref{cosinh'_a(0)}), we get:
\begin{equation*} \label{simplified_flex_equation}
u \tgh h \Big[ \cos \alpha \Big\{ \frac{\cosinh p \sinush q}{\tgh p}
+ \frac{\sinush p \cosinh q}{\tgh q} \Big\} - \frac{\sinush p
\cosinh q}{\tgh p} - \frac{\cosinh p \sinush q}{\tgh q} \Big] = 0.
\end{equation*}
Thus, the deformation under consideration of $\mathscr{S}$ is an
infinitesimal flex if and only if (\ref{v(u)}), (\ref{w(u)}) and
\begin{equation*} \label{prefinal_flex_equation}
\cos \alpha \Big\{ \frac{\cosinh p \sinush q}{\tgh p} +
\frac{\sinush p \cosinh q}{\tgh q} \Big\} = 2 \cosinh p \cosinh q
\end{equation*}
hold true. Hence, $\mathscr{S}$ allows the infinitesimal flex of the
form described in the beginning of this section if and only if $p$,
$q$, and $\alpha$ satisfy the following relation:
\begin{equation} \label{p_q_alpha}
\frac{\tgh p}{\tgh q} = \frac{1 \pm \sin \alpha}{\cos \alpha }.
\end{equation}

The so-constructed infinitesimal flex is nontrivial because the
distance between the poles $N(t)$ and $S(t)$ is not stationary.

\section{Calculating metric elements of
$\mathscr{T}(t)$} \label{paragraph_metric_elements}

Let's obtain formulae for the dihedral angles $\angle A(t)B(t)$,
$\angle N(t)A(t)$, and $\angle N(t)B(t)$ of the tetrahedron
$\mathscr{T}(t)$, which will be used in a proof of the theorem.

First we calculate the sines and cosines of the plane angles of
$\mathscr{T}(t)$.

Apply the Cosine Law for the Lobachevsky space to the triangle
$\triangle CA(t)N(t)$ to calculate the cosine of the angle $\beta
(t)$: $ \label{theorem_cos_beta(t)} \cosinh (h + tu) = \cosinh (p +
tv) \cosinh b(t) - \sinush (p + tv) \sinush b(t) \cos \beta (t) $.
Thus, taking into account~(\ref{cosinh_b(t)}) and formulae of
hyperbolic trigonometry, we get:
\begin{equation}
\label{cos_beta(t)} \cos \beta (t) = \frac{\sinush(p + tv) \cosinh
(h + tu)}{\sinush b(t)} = \frac{\sinush(p + tv) \cosinh (h +
tu)}{\sqrt{{\cosinh}^{2} (h + tu) {\cosinh}^{2} (p + tv) - 1}}.
\end{equation}
(Here and below $\sqrt{s}$ stands for a branch of the square root
that takes a positive real value for a positive real $s$.) To
calculate the sine of $\beta (t)$ we apply the Sine Law for the
Lobachevsky space \cite{Vinberg1988} to $\triangle CA(t)N(t)$:
\begin{equation*}
\label{theorem_sin_beta(t)} \frac{\sin \beta (t)}{\sinush (h + tu)}
= \frac{\sin \pi / 2}{\sinush b(t)} = \frac{1}{ \sqrt{{\cosinh}^{2}
(h + tu) {\cosinh}^{2} (p + tv) - 1}},
\end{equation*}
and therefore,
\begin{equation}
\label{sin_beta(t)} \sin \beta (t) = \frac{\sinush (h + tu)}{\sinush
b(t)} = \frac{\sinush (h + tu)}{ \sqrt{{\cosinh}^{2} (h + tu)
{\cosinh}^{2} (p + tv) - 1}}.
\end{equation}

Similarly, we obtain the formulae for the cosine and sine of the
angle $\varphi (t)$ in $\triangle CB(t)N(t)$:
\begin{equation}
\label{cos_phi(t)} \cos \varphi (t) = \frac{\sinush(q + tw) \cosinh
(h + tu)}{\sinush c(t)} = \frac{\sinush(q + tw) \cosinh (h +
tu)}{\sqrt{{\cosinh}^{2} (h + tu) {\cosinh}^{2} (q + tw) - 1}},
\end{equation}
\begin{equation}
\label{sin_phi(t)} \sin \varphi (t) = \frac{\sinush (h +
tu)}{\sinush c(t)} = \frac{\sinush (h + tu)}{ \sqrt{{\cosinh}^{2} (h
+ tu) {\cosinh}^{2} (q + tw) - 1}},
\end{equation}
for the cosine and sine of the angle $\mu (t)$ in $\triangle
CA(t)N(t)$:
\begin{equation}
\label{cos_mu(t)} \cos \mu (t) = \frac{\sinush(h + tu) \cosinh (p +
tv)}{\sinush b(t)} = \frac{\sinush(h + tu) \cosinh (p +
tv)}{\sqrt{{\cosinh}^{2} (h + tu) {\cosinh}^{2} (p + tv) - 1}},
\end{equation}
\begin{equation}
\label{sin_mu(t)} \sin \mu (t) = \frac{\sinush (p + tv)}{\sinush
b(t)} = \frac{\sinush (p + tv)}{ \sqrt{{\cosinh}^{2} (h + tu)
{\cosinh}^{2} (p + tv) - 1}},
\end{equation}
and for the cosine and sine of the angle $\nu (t)$ in $\triangle
CB(t)N(t)$:
\begin{equation}
\label{cos_nu(t)} \cos \nu (t) = \frac{\sinush(h + tu) \cosinh (q +
tw)}{\sinush c(t)} = \frac{\sinush(h + tu) \cosinh (q +
tw)}{\sqrt{{\cosinh}^{2} (h + tu) {\cosinh}^{2} (q + tw) - 1}},
\end{equation}
\begin{equation}
\label{sin_nu(t)} \sin \nu (t) = \frac{\sinush (q + tw)}{\sinush
c(t)} = \frac{\sinush (q + tw)}{ \sqrt{{\cosinh}^{2} (h + tu)
{\cosinh}^{2} (q + tw) - 1}}.
\end{equation}

The Cosine Law for the Lobachevsky space applied twice to the
triangle $\triangle A(t)CB(t)$ leads us to the formulae:
\begin{equation}
\label{cos_delta(t)} \cos \delta (t) = \frac{\cosinh(p + tv) \cosinh
a(t) - \cosinh(q + tw)}{\sinush(p + tv) \sinush a(t)} ,
\end{equation}
\begin{equation}
\label{cos_psi(t)} \cos \psi (t) = \frac{\cosinh(q + tw) \cosinh
a(t) - \cosinh(p + tv)}{\sinush(q + tw) \sinush a(t)}.
\end{equation}

From the Sine Law for the Lobachevsky space applied to $\triangle
A(t)CB(t)$, it follows that:
\begin{equation*}
\label{theorem_sin_ACB} \frac{ \sin \delta (t) }{\sinush(q + tw)} =
\frac{ \sin \alpha }{\sinush a(t)} = \frac{ \sin \psi (t) }{ \sinush
(p + tv) },
\end{equation*}
and thus the formulae
\begin{equation}
\label{sin_delta(t)} \sin \delta (t) = \frac{\sin \alpha \sinush (q
+ tw)}{\sinush a(t)} ,
\end{equation}
\begin{equation}
\label{sin_psi(t)} \sin \psi (t) = \frac{\sin \alpha \sinush (p +
tv)}{\sinush a(t)}
\end{equation}
hold true.

The Cosine Law for the Lobachevsky space three times applied to the
triangle $\triangle A(t)N(t)B(t)$ leads us to the formulae:
\begin{equation}
\label{cos_theta(t)} \cos \theta (t) = \frac{\cosinh a(t) \cosinh
c(t) - \cosinh b(t)}{\sinush a(t) \sinush c(t)} ,
\end{equation}
\begin{equation}
\label{cos_gamma(t)} \cos \gamma (t) = \frac{\cosinh a(t) \cosinh
b(t) - \cosinh c(t)}{\sinush a(t) \sinush b(t)},
\end{equation}
\begin{equation}
\label{cos_chi(t)} \cos \lambda (t) = \frac{\cosinh b(t) \cosinh
c(t) - \cosinh a(t)}{\sinush b(t) \sinush c(t)}.
\end{equation}

Taking into account (\ref{cosinh_b(t)})--(\ref{cosinh_a(t)}), we
calculate $\sinush a(t)$, $\sinush b(t)$, and $\sinush c(t)$ from
(\ref{cos_beta(t)})--(\ref{cos_chi(t)}):
\begin{equation*}
\label{sinush_a(t)} \sinush a(t) = \sqrt { {\cosinh}^{2} a(t) - 1} =
\sqrt { (\cosinh (p + tv) \cosinh (q + tw) - \sinush (p + tv)
\sinush (q + tw)\cos \alpha)^2 - 1 },
\end{equation*}
\begin{equation*}
\label{sinush_b(t)} \sinush b(t) = \sqrt { {\cosinh}^{2} b(t) - 1} =
\sqrt { (\cosinh (h + tu) \cosinh (p + tv) )^2 - 1 },
\end{equation*}
\begin{equation*}
\label{sinush_c(t)} \sinush c(t) = \sqrt { {\cosinh}^{2} c(t) - 1} =
\sqrt { (\cosinh (h + tu) \cosinh (q + tw) )^2 - 1 }.
\end{equation*}

The fact that the values of the angles in a hyperbolic triangle are
greater than $0$ and less than $\pi$ yields that the sines of the
angles of a hyperbolic triangle are nonnegative. Hence, $
\label{sin_theta(t)} \sin \theta (t) = \sqrt{1 - {\cos}^{2} \theta
(t) } $, $ \label{sin_gamma(t)} \sin \gamma (t) = \sqrt{1 -
{\cos}^{2} \gamma (t) }$, $\label{sin_chi(t)} \sin \lambda (t) =
\sqrt{1 - {\cos}^{2} \lambda (t) }$.

Consider the unit sphere $\Sigma$ centered at the vertex $A(t)$ of
$\mathscr{T}(t)$. Denote the points of the intersection of $\Sigma$
and the rays $\overrightarrow{A(t)C}$, $\overrightarrow{A(t)N(t)}$,
and $\overrightarrow{A(t)B(t)}$ by $C_A(t)$, $N_A(t)$, and $B_A(t)$
correspondingly. They determine a triangle $\triangle C_A(t) N_A(t)
B_A(t)$ which consists of the points of the intersection of $\Sigma$
and the rays emitted from $A(t)$ and passing through the points of
the face $\triangle CB(t)N(t)$ of $\mathscr{T}(t)$. By construction,
the angle of the spherical triangle $\triangle C_A(t) N_A(t) B_A(t)$
at the vertex $C_A(t)$ is equal to $\pi/2$, the angle at $N_A(t)$ is
equal to $\angle N(t)A(t)$, the angle at $B_A(t)$ is equal to
$\angle A(t)B(t)$, the length of the side $C_A(t) N_A(t)$ is equal
to $\beta (t)$, the length of $N_A(t) B_A(t)$ is equal to $\gamma
(t)$, and the length of $C_A(t) B_A(t)$ is equal to $\delta (t)$.

Similarly, we build a spherical triangle $\triangle C_B(t) N_B(t)
A_B(t)$. Its angle at the vertex $C_B(t)$ is equal to $\pi/2$, the
angle at $N_B(t)$ is equal to $\angle N(t)B(t)$, the angle at
$A_B(t)$ is equal to $\angle A(t)B(t)$, the length of the side
$C_B(t) N_B(t)$ is equal to $\varphi (t)$, the length of $N_B(t)
A_B(t)$ is equal to $\theta (t)$, and the length of $C_B(t) A_B(t)$
is equal to $\psi (t)$.

Applying the Cosine Law for the spherical space \cite{Vinberg1988}
twice to $\triangle C_A(t) N_A(t) B_A(t)$, we obtain the formulae:
\begin{equation*} \label{cos_AB(t)_CaNaBa} \cos \angle A(t)B(t)
= \frac{\cos \beta (t) - \cos \gamma(t) \cos \delta (t)}{\sin
\gamma(t) \sin \delta (t)},
\end{equation*}
\begin{equation*} \label{cos_NA(t)_CaNaBa} \cos \angle N(t)A(t)
= \frac{\cos \delta (t) - \cos \gamma(t) \cos \beta (t)}{\sin
\gamma(t) \sin \beta (t)}.
\end{equation*}

Again, applying the Cosine Law for the spherical space to $\triangle
C_B(t) N_B(t) A_B(t)$, we get:
\begin{equation*} \label{cos_NB(t)_CbNbAb} \cos \angle N(t)B(t)
= \frac{\cos \psi (t) - \cos \varphi(t) \cos \theta (t)}{\sin
\varphi(t) \sin \theta (t)}.
\end{equation*}

Now apply the Sine Law for the spherical space \cite{Vinberg1988} to
$\triangle C_A(t) N_A(t) B_A(t)$:
\begin{equation*} \label{theorem_of_sinuses_CaNaBa} \frac{\sin \angle N(t)A(t)}{\sin
\delta (t)} = \frac{\sin \angle A(t)B(t)}{\sin \beta (t)} =
\frac{\sin  \pi/2}{\sin \gamma (t)}.
\end{equation*}
Hence,
\begin{equation*} \label{sin_AB(t)_NA(t)_CaNaBa} \sin \angle A(t)B(t)
= \frac{\sin \beta (t)}{\sin \gamma (t)} \quad \mbox{and} \quad \sin
\angle N(t)A(t) = \frac{ \sin \delta (t)}{\sin \gamma (t)}.
\end{equation*}

Again, apply the Sine Law for the spherical space to $\triangle
C_B(t) N_B(t) A_B(t)$:
\begin{equation*} \label{theorem_of_sinuses_CnAnBn} \frac{\sin \angle N(t)A(t)}{\sin
\nu (t)} = \frac{\sin \angle N(t)B(t)}{\sin \mu (t)} = \frac{\sin
 \alpha }{\sin \lambda (t)}.
\end{equation*}
Thus,
\begin{equation*} \label{sin_NB(t)_CnAnBn} \sin \angle N(t)B(t)
= \sin \alpha \frac{ \sin \mu (t)}{\sin \lambda (t)}.
\end{equation*}

In the proof of the theorem given below we use also the following
three evident relations:
\begin{equation*} \label{NA(t)_prime} \frac{ d \angle N(t)A(t)}{d t}
= - \frac{\frac{ d }{d t} (\cos \angle N(t)A(t))}{\sin \angle
N(t)A(t)},
\end{equation*}
\begin{equation*} \label{NB(t)_prime} \frac{ d \angle N(t)B(t)}{d t}
= - \frac{\frac{ d }{d t} (\cos \angle N(t)B(t))}{\sin \angle
N(t)B(t)},
\end{equation*}
and
\begin{equation*} \label{AB(t)_prime} \frac{ d \angle A(t)B(t)}{d t}
= - \frac{\frac{ d }{d t} (\cos \angle A(t)B(t))}{\sin \angle
A(t)B(t)}.
\end{equation*}

\section{Proof of the theorem} \label{paragraph_proof_of_the_theorem}

Remind that, according to the Schl\"{a}fli formula for polyhedra in
the Lobachevsky 3-space \cite{Vinberg1988} of the curvature $-1$,
the equality
\begin{equation} \label{Schlafli_formula} dV = - \frac{1}{2}
\sum_{e} l_e d \theta_{e}
\end{equation}
holds true, where $dV$ stands for the variation of the volume of the
polyhedron, $l_e$ stands for the length of an edge $e$ of the
polyhedron, $d \theta_{e}$ stands for the variation of the dihedral
angle of the polyhedron attached to the edge $e$, and summation is
taken over all edges $e$ of the polyhedron.

Show that the polyhedron $\mathscr{S}(0)$ from the family of
suspensions $\mathscr{S}(t)$, $t \in (-1,1)$, constructed in
Section~\ref{paragraph_polyhedron_construction}, with parameters of
the tetrahedron $\mathscr{T}$
\begin{equation} \label{p_q_h_alpha}  p = \arctanh \frac{1}{2}, \,\,\,\,
q = \arctanh \frac{\sqrt{3}}{2}, \,\,\,\, h = \arctanh \frac{1}{2},
\,\,\,\, \alpha = \frac{\pi}{6} \,\,\,\, (\mbox{i. e.} \,\, n = 6)
\end{equation}
and the velocities of deformation
\begin{equation} \label{u_v_w}  u = \frac{\sqrt{3}}{4}, \,\,\,\,
v = - \frac{\sqrt{3}}{4}, \,\,\,\, w = - \frac{1}{4},
\end{equation}
can be taken as a polyhedron whose existence is asserted in the
theorem.

The suspension $\mathscr{S}(0)$ is not infinitesimally rigid because
$p$, $q$, and $\alpha$ from (\ref{p_q_h_alpha}) satisfy
(\ref{p_q_alpha}).

Let's verify that the nontrivial infinitesimal flex from
Section~\ref{paragraph_flex_condition} with the coefficients
(\ref{u_v_w}) can be taken as an infinitesimal flex whose existence
is stated in the theorem.

Using the Schl\"{a}fli formula (\ref{Schlafli_formula}) and taking
into account notations and remarks of
Section~\ref{paragraph_polyhedron_construction}, we see that the
variation of the volume of $\mathscr{S}(t)$ at $t = 0$ can be
written as follows:
\begin{equation} \label{Schlafli_formula_for_S(t)_part2}  dV_{\mathscr{S}(0)} = -
 12 \bigg{(}a(0) \frac{d \angle A(t)B(t)}{dt}(0) + b(0) \frac{d \angle N(t)A(t)}{dt}(0)
+ c(0) \frac{d \angle N(t)B(t)}{dt}(0) \bigg{)} dt.
\end{equation}
Substituting the values of parameters from (\ref{p_q_h_alpha}) and
(\ref{u_v_w}) into the formulae of
Sections~\ref{paragraph_flex_condition}
and~\ref{paragraph_metric_elements}, we sequentially find the
hyperbolic sines and cosines of the lengths of the edges and the
variations of the dihedral angles of the tetrahedron
$\mathscr{T}(t)$ at $t = 0$:
\begin{equation*}
\label{cosinh_a(t)_particular_case} \cosinh a(t) = \cosinh \bigg{(}
- \arctanh \frac{1}{2} \bigg{)} \cosinh \bigg{(} - \arctanh
\frac{\sqrt{3}}{2} \bigg{)} - \frac{\sqrt{3}}{2} \sinush \bigg{(} -
\arctanh \frac{1}{2} \bigg{)} \sinush \bigg{(} - \arctanh
\frac{\sqrt{3}}{2} \bigg{)},
\end{equation*}
\begin{equation*}
\label{cosinh_b(t)_c(t)_particular_case} \cosinh b(0) = \cosinh
\bigg{(} - \arctanh \frac{1}{2} \bigg{)} \cosinh \bigg{(} \arctanh
\frac{1}{2} \bigg{)}, \,\,\,\, \cosinh c(t) = \cosinh
\bigg{(}\arctanh \frac{1}{2} \bigg{)} \cosinh \bigg{(}- \arctanh
\frac{\sqrt{3}}{2} \bigg{)},
\end{equation*}
\begin{equation*}
\label{primes_of_dihedral_angles}  \frac{d \angle A(t)B(t)}{dt}(0) =
\frac{\sqrt{13}}{4}, \,\,\,\, \frac{d \angle N(t)A(t)}{dt}(0) =
\frac{\sqrt{7}}{4}, \,\,\,\, \frac{d \angle N(t)B(t)}{dt}(0) = -
\frac{\sqrt{13}}{4},
\end{equation*}
and thus, by (\ref{Schlafli_formula_for_S(t)_part2}),
\begin{equation*} \label{Schlafli_formula_for_S(t)_with_numbers_part1}
\frac{dV_{\mathscr{S}(0)}}{dt} =  -
 12 \bigg{[}\frac{\sqrt{13}}{4} \arch \Big{(} \cosinh
\bigg{(} - \arctanh \frac{1}{2} \bigg{)} \cosinh \bigg{(} - \arctanh
\frac{\sqrt{3}}{2} \bigg{)} -
\end{equation*}
\begin{equation*} \label{Schlafli_formula_for_S(t)_with_numbers_part23} -
\frac{\sqrt{3}}{2} \sinush \bigg{(} - \arctanh \frac{1}{2} \bigg{)}
\sinush \bigg{(} - \arctanh \frac{\sqrt{3}}{2} \bigg{)} \Big{)} +
\frac{\sqrt{7}}{4} \arch \Big{(} \cosinh \bigg{(} - \arctanh
\frac{1}{2} \bigg{)} \cosinh \bigg{(} \arctanh \frac{1}{2} \bigg{)}
\Big{)} -
\end{equation*}
\begin{equation*} \label{Schlafli_formula_for_S(t)_with_numbers_part4}
- \frac{\sqrt{13}}{4} \arch \Big{(} \cosinh \bigg{(} \arctanh
\frac{1}{2} \bigg{)} \cosinh \bigg{(} - \arctanh \frac{\sqrt{3}}{2}
\bigg{)} \Big{)} \bigg{]} =
\end{equation*}
\begin{equation*} \label{Schlafli_formula_for_S(t)_with_numbers_part567}
- 3 \Big{[} \sqrt{7} \arch \frac{4}{3} + \sqrt{13} \Big{(} \arch
\frac{5}{2 \sqrt{3}} - \arch \frac{4}{\sqrt{3}}\Big{)} \Big{]}  = -
3  \Big{[} \sqrt{7} \ln \frac{4 + \sqrt{7}}{3} + \sqrt{13} \ln
\frac{7 - \sqrt{13}}{6} \Big{]} <
\end{equation*}
\begin{equation*} \label{Schlafli_formula_for_S(t)_with_numbers_part8}
< - \frac{3 \sqrt{7}}{8} \Big{[} 8 \ln \frac{4 + \sqrt{7}}{3} + 11
\ln \frac{7 - \sqrt{13}}{6} \Big{]} = - \frac{3 \sqrt{7}}{8} \ln
\Big{[} \Big{(} \frac{4 + \sqrt{7}}{3} \Big{)}^8 \Big{(} \frac{7 -
\sqrt{13}}{6} \Big{)}^{11} \Big{]} < 0. \quad \square
\end{equation*}

\section{Concluding remarks} \label{paragraph_final_conclusions}

Using notations of Section~\ref{paragraph_proof_of_the_theorem}, we
determine the integral mean curvature of a polyhedron
$\mathscr{S}(t)$ in the 3-space as follows:
\begin{equation*} \label{Mean_curvature_formula} M(\mathscr{S}(t)) = \frac{1}{2}
\sum_{e} l_e(t) ( \pi - \theta_{e}(t)).
\end{equation*}

R.~Alexander \cite{Alexander1985} proved that the integral mean
curvature of any polyhedron in the Euclidean 3-space is stationary
under every its infinitesimal flex.

The lengths of the edges of the suspension $\mathscr{S}(t)$ are
stationary under the infinitesimal flex of $\mathscr{S}(t)$ from
Section~\ref{paragraph_flex_condition}. Hence, the variation of the
integral mean curvature of $\mathscr{S}(t)$ at $t = 0$ is equal to
the variation of the volume $dV_{\mathscr{S}(0)}$. Therefore, the
proof of our theorem automatically implies that the variation of the
integral mean curvature for the infinitesimal flex of
$\mathscr{S}(t)$ constructed above is not equal to zero. Thus, the
integral mean curvature of an infinitesimally nonrigid polyhedron is
not always stationary in the Lobachevsky space as well as in the
spherical space but is always stationary in the Euclidean space.

The author is thankful to Victor Alexandrov for his help in the
preparation of the present paper and, generally, for his support of
the author's studies and research in mathematics.

Dmitriy Slutskiy

Sobolev Institute of Mathematics of the SB RAS,

4 Acad. Koptyug avenue, 630090 Novosibirsk, Russia

and

Novosibirsk State University,

2 Pirogova Street, 630090, Novosibirsk, Russia

slutski@ngs.ru

\end{document}